\documentclass{article}
              \usepackage{amsfonts}

\newtheorem{Th}{Theorem}

\newtheorem{Lemma}{Lemma}
\newtheorem{Cor}{Corollary}
\newtheorem{Def}{Definition}
\newtheorem{Rem}{Remark}
\newtheorem{Con}{Conjecture}

\newcommand{\eqdef}{\stackrel{{\rm def}}{=}}

\newcommand{\diag}{\mbox{\rm diag}}

\newcommand{\Id}{\mbox{\rm Id}}

\def\diag{\mbox{\rm diag}}

\title{On projectively equivalent metrics near points of bifurcation}
\author{Vladimir S. Matveev\thanks{Mathematisches Institut, Universit\"at 
Freiburg, 79104 Germany \newline
 matveev@minet.uni-jena.de}}
\date{}
\begin{document}
\maketitle

\noindent{\bf Keywords: }{ \ 
\small  Projectively Equivalent Metrics, Geodesically Equivalent Metrics,
 Integrable Systems, Levi-Civita Coordinates}

\vspace{2ex}

\noindent{\bf MSC2000: }{ \ \small 37J35,  70H06,    53D25,   53B10
 53C24,  53C15,  37J30,  53A20}

\begin{abstract}
 Let Riemannian metrics $g$ and $\bar g$ on a  connected manifold 
  $M^n$ 
 have the same geodesics (considered as unparameterized curves). 
 Suppose the eigenvalues of one metric with 
 respect to the other are all different at  a  point.
 Then, by the famous Levi-Civita's Theorem, the metrics have a certain 
standard form near the point. 
 
   Our main result  is a  
   generalization  of  Levi-Civita's Theorem  for 
    the points where   the eigenvalues of  one metric with 
 respect to the other  bifurcate.   
\end{abstract}

\section{Introduction}
\subsection{ Metrics with the same geodesics}
\begin{Def}\hspace{-2mm}{\bf .} {\rm 
Two Riemannian metrics $g$ and   $\bar g$  on $M^n$  
  are called 
 {\em projectively  equivalent}, 
 if they have the same geodesics 
considered as unparameterized curves.}
\end{Def}

Trivial examples of  projectively equivalent metrics can be obtained by 
considering proportional metrics 
 $g$ and $C\cdot  g$, where $C$ is a positive 
constant.

\begin{Def}\hspace{-2mm}{\bf .} {\rm 
 The Riemannian metrics $g$ and $\bar g$ 
are said to be {\em strictly non-proportional}  at $x\in M^n$, if 
the eigenvalues of $g$ with 
respect to $\bar g$  are all different at $x$. }
\end{Def}

In the present paper we study projectively  equivalent metrics which 
are strictly non-proportional at least at one point of the manifold.  
It is a very classical material. 
 In 1865, Beltrami \cite{Beltrami} found the first examples 
 and  formulated a problem of finding all pairs
   of projectively equivalent metrics. 
Locally the problem almost has been solved by 
Dini \cite{Dini} for dimension two and Levi-Civita  \cite{Levi-Civita} 
for  an arbitrary dimension: they  
 obtained a local description of projectively equivalent metrics near 
the points where the  eigenvalues of one metric with 
respect to the other do not bifurcate: 
\begin{Th}[Levi-Civita \cite{Levi-Civita}]\hspace{-2mm}{\bf .} \label{Levi-Civita} Consider two  Riemannian metrics 
 on an open subset $U^n\subset M^n$. 
 Suppose the metrics are strictly-non-proportional at every  point $x\in U^n$. 

Then the metrics are projectively equivalent on $U^n$,  
if and only 
if for every   point  $x\in U^n$ 
there exist coordinates $x_1, x_2, ... , x_n$ in some  neighborhood 
of   $x$ 
such that in these coordinates the metrics have the 
following model  form: 
\begin{eqnarray}
ds_{g_{model}}^2& = & \
\Pi_1dx_1^2+
\ \Pi_2dx_2^2\ +\cdots+
\ \ \Pi_{n}dx_{n}^2,\label{Canon_g}\\
ds_{\bar g_{model}}^2&=&\rho_1\Pi_1dx_1^2+
\rho_2\Pi_2dx_2^2+\cdots+
\rho_{n}\Pi_{n}dx_{n}^2 \label{Canon_bg},
\end{eqnarray}
where the functions 
$\Pi_i $ and $\rho_i$
 are given  by  
\begin{eqnarray} 
\Pi_i & \eqdef & \label{Pi}
(\lambda_i-\lambda_1)(\lambda_i-\lambda_2)
\cdots(\lambda_{i}-\lambda_{i-1})(\lambda_{i+1}-
\lambda_i)\cdots
(\lambda_{n}-\lambda_i), \\ 
 \rho_i  & \eqdef &
 \frac{1}{\lambda_1\lambda_2...\lambda_{n-1}}\frac{1}{\lambda_i}.  \nonumber
\end{eqnarray}
where,    for any $i$, the function  
$
\lambda_i
$ 
is a smooth function of the variable $x_i$.  
\end{Th}

However, the global behavior of 
projectively equivalent metrics
is not understood completely.

 One of the main 
 difficulties for global description of projectively equivalent 
 metrics was  the absence  of local 
  description of  projectively equivalent metrics near the 
 points where the eigenvalues of  one metric with 
   respect to the other  bifurcate.  As it has been proven in \cite{commute}, 
   if the manifold is not covered by the torus, it must have such points.

Our main result is  a 
   description of (strictly non-proportional at least at one point) 
  projectively equivalent metrics near  points
   where the eigenvalues of  one metric with 
   respect to the other  bifurcate.  
 The description needs some preliminary work and will be formulated 
  at the very end of  the  paper,  see Theorem~\ref{local} in 
Section~\ref{2.6}.  Here, 
  we would like to note that the points of bifurcation behave quite 
  regularly: the multiplicity of an eigenvalue 
is at most three; the points where the multiplicity of an eigenvalue is at 
least two are organized in totally-geodesic submanifold $U^{n-2}$ 
of co-dimension two; the points where the multiplicity of the 
 eigenvalue is three are organized in a closed totally-geodesic hypersurface
 (submanifold of co-dimension one) of $U^{n-2}$.

This local description  is the main tool  of 
 the following 

\begin{Th}[Topalov, Matveev]\hspace{-2mm}{\bf .}\label{main}
Suppose $M^n$ is a connected closed manifold. 
Let Riemannian metrics  
$g$  and $\bar g$ on $M^n$ 
 be projectively 
equivalent. Suppose they are  strictly non-proportional 
at least at 
one point. Then the manifold can be finitely covered by a  product of 
 spheres. 
\end{Th}

We will neither prove nor explain this theorem in this paper. The 
proof will appear elsewhere, in a joint paper with 
 P. Topalov. 

As we will show in 
 Section~\ref{2.3}, see Corollary~\ref{product}, a 
 product of spheres always admits 
   projective equivalent metrics strictly non-proportional 
at least at 
one point.

\subsection{ Geodesic rigidity conjecture}

\begin{Con}\hspace{-2mm}{\bf .}\label{conjecture}
Suppose $M^n$ is closed connected. Let   Riemannian metrics
$g$ and $\bar g$ on $M^n$ be  
 projectively equivalent and be nonproportional. 
Then the manifold can be covered by the sphere,  
or it admits a local product structure. 
\end{Con} 

Definition  of  the notion 
``local product structure''  is in Section~\ref{2.2}.    Roughly speaking, 
local product structure is a  metric and two orthogonal foliations of 
complimentary dimensions such that locally all three objects 
 look as they came from the  
direct product of two Riemannian manifold. 

Conjecture~\ref{conjecture} is true for dimensions two  \cite{surface,ERA,dedicata}
and three \cite{short,threemanifolds}. Theorem~\ref{main} shows 
that  conjecture is true, if we assume in addition that 
the metrics are strictly-non-proportional at least at one point. 
The following theorem shows that the conjecture is true on the level of 
fundamental groups: 
\begin{Th}[\cite{hyperbolic}]\hspace{-2mm}{\bf .} 
Let $M^n$ be a closed connected   manifold. Suppose  two
 nonproportional  Riemannian metrics $g$, $\bar g$ on $M^n$  
are   projectively equivalent.
If  the  fundamental group   of $M^n$ is infinite, then 
  there exist  a local product structure on $M$. 
\end{Th}

If a manifold admits a local product structure, it admits two non-proportional 
projectively equivalent metrics. The direct 
  product   of manifolds admitting strictly-non-proportional projectively equivalent metrics also admit strictly-non-proportional projectively equivalent metrics (see Section~\ref{2.3} and Lemma~\ref{inverse} there).  Every sphere admits strictly-non-proportional projectively equivalent metrics (essentially
\cite{Beltrami}; see also  
Corollary~\ref{product} in Section~\ref{2.3}).  

From the 
other side,  not every manifold 
 covered by the  sphere can carry
 non-proportional   projectively equivalent metrics.
 The lowest dimension where 
 such phenomena is possible  is three: 
\begin{Th}[\cite{short},\cite{threemanifolds}]\hspace{-2mm}{\bf .}\label{dim3}
 A closed  connected 3-manifold admits   non-proportional projectively equivalent 
metrics if and only if it is homeomorphic to a lens space or to a Seifert manifold 
with zero Euler number.
\end{Th}
 Then, the Poincare homology sphere 
(which is certainly covered by $S^3$ \cite{Kirbi}) does not 
 admit   non-proportional projectively equivalent metrics.

\subsection{ Integrable systems in the theory of projectively equivalent
 metrics}

 There exist   two classical  techniques  for working with projectively
  equivalent metrics:

One is  due to Dini \cite{Dini}, Lie \cite{Lie}, 
 Levi-Civita \cite{Levi-Civita}, Weyl \cite{Weyl}  and Eisenhart \cite{Eisenhart}, and was 
actively developed   by Russian and Japan geometry schools in the 50th-70th.  
Projective  equivalence of two metrics is equivalent to a certain 
system of quasi-linear
 differential equation on the entries of these metrics. Combining 
 this system with additional geometric assumptions  (mostly written as a tensor formula) 
one can deduce topological restrictions on the manifold, 
see the surveys \cite{Mikes2}.

The second classical technique is due to   E. Cartan \cite{cartan}  
 and was very actively 
developed by  the French  geometry  school. The main observation is that, 
 on the level of connections, projective  equivalence is a very easy condition. 
In particular,  given a connection, it is easy  to describe 
 all projectively equivalent connections. So the goal is to 
understand which of these connections are Levi-Civita connections of a 
Riemannian  metric.  

Unfortunately, both techniques are local, and all global (when the 
manifold is closed) result obtained with the help of these 
techniques require additional strong geometric assumptions
(i.e.,    the metric  is assumed to 
 be K\"ahler or semisymmetric or generally semisymmetric or T-generalized semisymmetric or Einstein or 
 of constant curvature   or Ricci-flat or recurrent or admitting a
 concircular vector field or admitting 
 a torse-forming vector field).

New methods  for global 
investigation of projectively 
 equivalent metrics  have been suggested in 
\cite{MT,TM}. The main observation of  \cite{MT,TM} is that the existence of 
$\bar g$ projectively equivalent to $g$ allows one to construct commuting 
integrals for the geodesic flow of $g$, see 
 Theorem~\ref{integrability} in Section~\ref{2.1}.

When the metrics are strictly non-proportional at least at one point, 
 the number of functionally 
independent integrals equals  the dimension of the manifold so that 
 the 
geodesic flow of $g$ is 
 Liouville-integrable  and we can apply 
the well-developed machinery of integrable geodesic flows. 
In particular, in dimension two, Theorem~\ref{main} and 
Conjecture~\ref{conjecture} follow directly from the description 
of metrics on surfaces 
with quadratically integrable geodesic flows from \cite{BMF,BF,Kio,Kol}. 

For arbitrary dimension, the  
following theorem has been obtained in \cite{added_reference} (the three-dimensional version  is due to \cite{starrheit,japan}) 
by combining  the ideas  applied  in  \cite{Taimanov} for analytically-integrable geodesic flows   
with  technique developed in  \cite{memo} for quadratically-integrable
geodesic flows:
\begin{Th}[\cite{added_reference}]\hspace{-2mm}{\bf .}\label{dubrovin}
Let  $M^n$ be  a connected closed manifold and  
$g$,  $\bar g$ 
 be projectively 
equivalent Riemannian metrics on $M^n$.  Suppose 
there exists a point of the manifold where 
the number of different  eigenvalues of $g$ with respect to $\bar g$
 equals $n$. Then the following holds:
\begin{enumerate} 
\item 
 The first Betti number $b_1(M^n)$ is 
not  greater than $n$. 

\item The fundamental group of the manifold is virtually Abelian.

\item If   there exists 
a point where 
the number of different  eigenvalues of $g$ with respect to $\bar g$
is less  than $n$,  
then $b_1(M^n)<n$.

\item If  there exists  no such point, then  
       $M^n$ can be covered by  the torus $T^n$. 
\end{enumerate}
\end{Th}

\noindent For dimensions three, 
 Theorem~\ref{main} follows  immediately from 
Theorem~\ref{dubrovin} modulo the  Poincar\'e conjecture.

More precisely, by the second statement of Theorem~\ref{dubrovin}, 
 there exists a finite cover with
 the  fundamental   group isomorphic to 
$\underbrace{Z\times Z\times ...\times Z}_k.$  
By \cite{Epstein},  the fundamental 
group of a compact 3-manifold  has no  subgroup isomorphic to 
$Z\times Z\times Z\times Z$,  so  $k$ is as most $3$. 

By  a  result  of Reidemeister (which has been formulated and proven 
in    \cite{Epstein1} and \cite{Thomas}), 
 it follows that  
$Z\times Z$ can not be the  fundamental group of a 
closed 3-manifold. 
Thus,  $k$ is  either 3 or 1 or  0.
  If $k=3$, then, by the fourth 
 statement of Theorem~\ref{dubrovin}, $M^3$ is covered by the torus $S^1\times S^1\times S^1$.    
If $k=0$, the fundamental group of $M^3$ is finite so that  
(modulo the Poincar\'e conjecture) it is covered
 by  $S^3$. 
If $k=1$, then,  modulo the Poincar\'e conjecture and in view of results
 of \cite{Whitehead}, $M^3$ is covered by the 
 product  $S^2\times S^1$, see \cite{japan}  for details.

\subsection{Acknowledgements}
I  would like to thank 
 W. Ballmann, V. Bangert, A. Bolsinov, K. Burns, N. A'Campo,  
     A. Fomenko, M. Gromov, U. Hamenst\"adt, K. Kiyohara, V. Lychagin, 
 A. Naveira,   S. Matveev,  P. Topalov,  and K. Voss  
 for  fruitful discussions.   
 I am very  grateful to  The Max-Planck Institute for Mathematics (Bonn), 
 The Institut des Hautes \`Etudes Scientifiques,   
   and   The Isaac Newton Institute for Mathematical Sciences
 for their hospitality and The European Post-Doctoral Institute 
 for   partial financial support.  My research 
 at INIMS has been  supported   by EPSRC  grant GRK99015.

\section{ Local description near bifurcations } \label{1}
Within this section, we assume  that Riemannian metrics $g$ and $\bar g$ 
are projectively   equivalent, and that the manifold is connected.  
Our goal  is to give  a local description of (strictly-non-proportional at
 least at one point) projectively equivalent metrics $g$, $\bar g$  
near the points where the eigenvalues of one  metric with respect to
 the other bifurcate (i.e. are not  all different).   Sections  
\ref{2.1},\ref{2.2},\ref{2.3} shows that it is sufficient to do these for 
two- and three-dimensional manifolds only; we do this in 
 Sections    \ref{2.4},\ref{2.5}  (Theorems~\ref{local2},\ref{local3}) 
 and combine all results in Section~\ref{2.6}, see Theorem~\ref{local}.

\subsection{The multiplicity of an eigenvalue is not  greater than three.}
 \label{2.1}

\noindent 
Let $g=(g_{ij})$ and  $\bar g=(\bar g_{ij})$ be Riemannian metrics on 
a manifold $M^n$. 
Consider the (1,1)-tensor  $L$ 
given by the formula   
\begin{eqnarray}
L^i_j &\eqdef & \left(\frac{\det(\bar g)}{\det(g)} \right)^\frac{1}{n+1} \bar g^{i\alpha} g_{\alpha j}. \label{l}
\end{eqnarray}
Then, $L$ determines the family $S_t$, $t\in R$, of $(1,1)$-tensors 
\begin{equation}\label{st}
 S_t\eqdef \det(L - t\ \mbox{\rm Id})\left(L-t\ \mbox{\rm Id}\right)^{-1}. 
 \end{equation}
\begin{Rem}\hspace{-2mm}{\bf .} 
Although $\left(L-t\ \Id\right)^{-1}$ is not defined for 
$t$ 
lying 
in the spectrum of $L$, the tensor  $S_t$  is well-defined 
for every  $t$.  Moreover,   $S_t$
is a  polynomial
 in $t$ of degree $n-1$ 
with coefficients being  (1,1)-tensors.  
\end{Rem}
\noindent We will identify the tangent and  cotangent bundles of $M^n$ by $g$. 
This identification allows us to transfer the 
 natural  Poisson structure from  $T^*M^n$ to   $TM^n$.

\begin{Th}[\cite{MT,TM}]\hspace{-2mm}{\bf .}
\label{integrability}
 If $g$, $\bar g$ are projectively 
 equivalent,   
then, for every  $t_1,t_2\in R$, the functions 
\begin{equation}\label{integral}
I_{t_i}:TM^n\to R, \ \ I_{t_i}(v)\eqdef g(S_{t_i}(v),v)
\end{equation}
are commuting integrals for the geodesic flow
 of  $g$. 
\end{Th}

 Since $L$ is  self-adjoint  with 
respect to $\bar g$,  the  eigenvalues of $L$  are real.  
At every point $x\in M^n$, let  
 us denote by  $\lambda_1(x)\le ... \le \lambda_n(x)$  the eigenvalues 
 of $L$ at the point.

\begin{Rem}\hspace{-2mm}{\bf .}
The notation $\lambda$ for the eigenvalues of $L$ is compatible with 
the notations used inside Levi-Civita's Theorem:  
in Levi-Civita coordinates from Theorem~\ref{Levi-Civita}, the tensor 
$L$  is given by the diagonal matrix  
\begin{equation}\label{lambda}
\diag(\lambda_1,\lambda_2,...,\lambda_n).  
\end{equation}
\end{Rem}

\begin{Def}\hspace{-2mm}{\bf .}
A Riemannian manifold is called {\bf geodesic}, if every  two points 
can be connected by a geodesic
\end{Def} 
 
\begin{Cor}[\cite{dedicata,hyperbolic}]\hspace{-2mm}{\bf .}\label{ordered12}
Let $(M^n, g)$ be a geodesic  Riemannian manifold. 
Let a Riemannian metric $\bar g$ on $M^n$ be projectively equivalent to $g$. 
 Then,   for every  $i\in \{1,... ,n-1\} $, for every  $x,y\in M^n$, the
 following holds:
\begin{enumerate}
\item $\lambda_i(x)\le \lambda_{i+1}(y)$. 

\item  If $\lambda_i(x)< \lambda_{i+1}(x)$, 
then $\lambda_i(z)< \lambda_{i+1}(z)$
for almost every point $z\in M^n$.
\end{enumerate}
\end{Cor}

\noindent In order to prove  Corollary~\ref{ordered12},  
we need the following technical lemma. 
For every fixed $v=(\xi_1,\xi_2, ... ,\xi_n)\in T_xM^n$,
 the function (\ref{integral}) is 
 a  polynomial  in $t$. 
Consider the roots of this  polynomial. 
 From the proof of Lemma~\ref{technical}, 
it will be clear that they are real. We denote them by 
 $$
t_1(x,v)\le t_2(x,v)\le ... \le  t_{n-1}(x,v).
$$ 
 \begin{Lemma}[\cite{dedicata,hyperbolic}]\hspace{-2mm}{\bf .} \label{technical} The following holds
 for every  $i\in \{1,...,n-1\}$: 
\begin{enumerate}
 \item For every  $v \in T_x M^n$, 
  $$
    \lambda_i(x)\le t_i(x,v) \le \lambda_{i+1}(x).   
  $$
     In particular,  if $\lambda_i(x) = \lambda_{i+1}(x)$, then 
$t_i(x,v)=\lambda_i(x) = \lambda_{i+1}(x)$. 

  \item If $\lambda_i(x) < \lambda_{i+1}(x)$, then for every   $\tau\in R $ 
the 
   Lebesgue 
measure of the set
 $$
V_\tau\subset   T_x M^n, \ \ V_\tau\eqdef \{ v\in T_xM^n: t_i(x,v) = \tau\},
$$
  is zero. 
\end{enumerate}
 \end{Lemma}

\noindent{\bf  Proof 
of Lemma~\ref{technical}:} 
By definition, the tensor $L$ is self-adjoint with respect to $g$. 
Then, for every  $x\in M^n$, 
 there exists  "diagonal" coordinates  in $T_xM^n$ where  the 
 metric $g$ is given by the    diagonal matrix 
 $\diag(1,1, ... ,1)$  and 
the tensor  $L$ is given by the 
  diagonal  matrix $\diag(\lambda_1,\lambda_2, ... ,\lambda_n)$. 
Then, the tensor 
 (\ref{st}) reads:
\begin{eqnarray*}
S_t& = & \det(L-t\Id)(L-t\Id)^{(-1)} \\
&=& \diag(\Pi_1(t), \Pi_2(t),
... ,\Pi_n(t)), 
\end{eqnarray*}
where the polynomials  $\Pi_i(t)$ are  given by the formula 
$$
\Pi_i(t)\eqdef (\lambda_1-t)(\lambda_2-t) ... (\lambda_{i-1}-t)(\lambda_{i+1}-t) ... (\lambda_{n-1}-t)(\lambda_n-t).   
$$
Hence, for every  $v=(\xi_1,...,\xi_n)\in T_xM^n$, 
the polynomial $I_t(x,v)$ is given by 
\begin{equation}\label{33}
I_t=\xi_1^2\Pi_1(t)+ \xi_2^2\Pi_2(t)+
... +\xi_n^2\Pi_n(t).
\end{equation} 
Evidently,  the coefficients of the polynomial $I_t$
depend continuously 
on the eigenvalues $\lambda_i$ and on the components 
 $\xi_i$. 
Then, it is sufficient to prove the first statement of the lemma assuming 
that 
the eigenvalues  $\lambda_i$ are all different and that  $\xi_i$
 are 
non-zero.  For every  $\alpha\ne i$, 
 we evidently 
have $\Pi_\alpha(\lambda_i)\equiv 0$. 
Then, 
$$
I_{\lambda_i} =\sum_{\alpha=1}^{n}\Pi_\alpha (\lambda_i) \xi_\alpha^2 =
\Pi_i(\lambda_i) \xi_i^2. 
$$
Hence  $I_{\lambda_i}(x,v)$ 
and  $I_{\lambda_{i+1}}(x,v)$ have different signs.
 Hence,
 the open interval $]\lambda_i,\lambda_{i+1}[$ contains  
a root of the polynomial $I_t(x,v)$. The degree of the polynomial $I_t$ is equal 
$n-1$; we have
 $n-1$ disjoint intervals; each of these
 intervals contains at least one 
root so that all roots are real and the  $i$th  root 
 lies between 
$\lambda_i$ 
and $\lambda_{i+1}$. 
  The first statement of the lemma is proved. 

Let us prove the second 
statement of Lemma~\ref{technical}.   
Suppose $\lambda_i<\lambda_{i+1}$.
Let first $\lambda_i<\tau< \lambda_{i+1}$. 
Then, the set 
$$  V_\tau\eqdef \{ v\in T_xM^n: t_i(x,v) = \tau\},$$
consists of the points $v$ where the function 
$I_{\tau}(x,v)\eqdef (I_t(x,v))_{|t= \tau}$ 
is zero; then   it is 
  a nontrivial quadric in $T_xM^n\equiv R^n$ and its measure is
 zero.

Let $\tau$ be 
one of the endpoints of the interval
$[\lambda_i, \lambda_{i+1}]$. 
Without loss of generality, we can suppose $\tau = \lambda_i$. 
Let $k$ be the multiplicity of the eigenvalue $\lambda_i$. Then, every  
coefficient $\Pi_\alpha(t)$ of the quadratic form 
 (\ref{33}) has  the  factor 
$(\lambda_i - t)^{k-1}$. Hence, 
$$
\hat I_t \eqdef  \frac{I_t}{(\lambda_i - t)^{k-1}}
$$
is a polynomial in $t$ and $\hat I_\tau$ is a nontrivial quadratic form.   
Evidently, for every  point   $v \in V_\tau$, we have 
$\hat I_\tau(v) = 0$ 
so that the set 
$V_\tau$ is a subset of 
a nontrivial quadric in  
 $T_xM^n$ and  its measure  is zero.  
  Lemma~\ref{technical}   is proved.

{\noindent \bf Proof of Corollary~\ref{ordered12}:}
The first statement of Corollary~\ref{ordered12} follows immediately from the 
first statement of Lemma~\ref{technical}:  
Let us join the points $x,y\in M^n$  by a geodesic
$\gamma:R\to M^n$, $\gamma(0)=x, \ \gamma(1)=y$. 
Consider the one-parametric family of integrals
$ 
I_t(x,v)
$ 
and the  roots 
$$
t_1(x,v) \le t_2(x,v) \le ... \le t_{n-1}(x,v).
$$

 By Theorem~\ref{integrability}, 
 each  root $t_i$ is   constant  
on every  orbit $(\gamma,\dot\gamma)$ of 
the geodesic flow of $g$ so that 
$$ 
t_i(\gamma(0),\dot\gamma(0)) =  t_i(\gamma(1),\dot\gamma(1)).
$$
Using Lemma~\ref{technical}, we obtain 
$$
\lambda_i(\gamma(0))\le t_i(\gamma(0),\dot\gamma(0)), \ \ \ \mbox{and}\ \ \ \
 t_i(\gamma(1),\dot\gamma(1)) \le \lambda_{i+1}(\gamma(1)).
$$
Thus  
$\lambda_i(\gamma(0))\le \lambda_{i+1}(\gamma(1))$ and the first statement 
of Corollary~\ref{ordered12} is proved.

Let us prove the second statement of Corollary~\ref{ordered12}. 
 Suppose 
 $\lambda_i(y) = \lambda_{i+1}(y)$ 
 for every point $y$ of  some  subset 
$V\subset M^n$. 
Then,   
the value of $\lambda_i$ is a  
constant 
(independent
 of $y\in V$). Indeed, by the first statement
 of Corollary~\ref{ordered12}, 
$$
\lambda_i(y_0)\le \lambda_{i+1}(y_1)   \  \ \mbox{and} \  \
\lambda_i(y_1) \le  \lambda_{i+1}(y_0), 
$$
so that $
\lambda_i(y_0)=\lambda_i(y_1)=\lambda_{i+1}(y_1)=\lambda_{i+1}(y_0)$ 
for every  $y_0, y_1\in V$. 

We denote this constant by $\tau$. 
 Let us join the point $x$ 
with every point 
of 
$V$ by all possible geodesics. Consider the set  
$V_\tau \subset T_xM^n$ 
of the 
initial velocity vectors (at the point $x$)
 of these geodesics.

 By the first statement of   Lemma~\ref{technical},
 for every  geodesic $\gamma$ passing
 through at least one  point 
of $V$,  the value  
$t_i(\gamma,\dot\gamma)$
 is equal to $\tau$. 
By the second  statement of   Lemma~\ref{technical}, the measure of the 
set $V_\tau $ is zero. Since the set 
$V$ lies in the image of the exponential mapping of the set $V_\tau$, 
 the measure of the set $V$ is also zero.
Corollary~\ref{ordered12}  is proved.

\begin{Cor}\hspace{-2mm}{\bf .}\label{multiplicity}
 Suppose $M^n$ connected. 
Let $g$ and $\bar g$ be projectively equivalent on $M^n$. Suppose they are
 strictly non-proportional at least at one point. Then the following holds:
\begin{enumerate} 
\item The metrics are strictly non-proportionally at almost every point of $M^n$. 

\item At every point of $M^n$,  the multiplicity of 
 any eigenvalue of $L$ is at most three. 
\end{enumerate}
\end{Cor} 

{\noindent} {\bf Proof:} Suppose the metrics are simply non-proportional at 
 $x_1$. Since the manifold is connected, for every  point $x\in M^n$ 
there exists a sequence of open convex balls $B_i\subset M^n$ such 
 that the  first ball $B_1$ contains the point $x_1$, the last ball $B_m$ 
 contains the point $X$ and the intersection  $B_i\cap B_{i+1}$, $i<m$, 
 is not empty. By the second  statement of Corollary~\ref{ordered12}, 
 the metrics are strictly non-proportional at 
 almost every  points  of  $B_1$. Then there are strictly non-proportional at 
 almost every  points  of $B_1\cap B_{2}$. Then, there exists a point of 
$B_2$ where the metrics are strictly non-proportional. 
Hence, by the second  statement of Corollary~\ref{ordered12},
 they   are strictly non-proportional at   almost every  points  of  $B_2$. Iterating this argumentation $m-2$ times,
  we   obtain that the metrics    are strictly non-proportional at 
 almost every  points  of  $B_m$.    Since the point $x$ is arbitrary, 
 the metrics are  strictly non-proportional at 
 almost every  points  of  $M^n$.

 Then, the multiplicity of an eigenvalue at the point $x$ is not 
greater than $3$: indeed, if we suppose that $\lambda_k(x)=\lambda_{k+3}(x)$, 
then, by the first statement  of Corollary~\ref{ordered12},
$\lambda_{k+1}=\lambda_{k+2}$ at each point of 
$B_m$, which contradicts that the metrics are   strictly non-proportional at 
 almost every  points  of  $B_m$. Corollary~\ref{multiplicity}  is proved.

\subsection{Splitting procedure} \label{2.2}
\noindent The goal of  the next  two sections is to show that for a local description of 
  projectively equivalent metrics strictly-non-proportional at least at 
 one point 
 it is sufficient to describe them on two- and three-dimensional manifold only.
We will need the following statement. 

\begin{Cor}[\cite{Benenti}]\hspace{-2mm}{\bf .} \label{nijenhuis}
Suppose the Riemannian metrics $g, \bar g$  on $M^n$ 
are projectively equivalent. 
Then  the Nijenhuis torsion of  $L$ vanishes.  
\end{Cor}

A self-contained proof of Corollary~\ref{nijenhuis} can be found in 
\cite{Benenti}; here we prove the theorem 
  assuming that the manifold is connected and 
the  metrics are strictly-non-proportional at least at one point, which is 
sufficient for our goals.

\noindent {\bf Proof:} Nijenhuis torsion is a tensor, 
 so it is sufficient to check its vanishing at almost every point. 
 By Corollary~\ref{ordered12}, almost every point of $M^n$ is stable. 
In the  Levi-Civita coordinates from Theorem~\ref{Levi-Civita}, the tensor 
$L$  is given by the diagonal matrix  
$$
\diag(\lambda_1,...,,\lambda_n). 
$$
Since the eigenvalue $\lambda_i$  depends  on  the variable  
$x_i$ only, the  Nijenhuis torsion of $L$ is  zero \cite{Haantjes}. 
Corollary~\ref{nijenhuis}
 is proved. 

\begin{Def}\hspace{-2mm}{\bf .} \label{lp} 
A local-product structure on $M^n$ is the triple $(h, B_r, B_{n-r})$,  
where  $h$ is a Riemannian metrics and $B_r$, 
$B_{n-r}$ are transversal foliations of dimensions   $r$ and $n-r$, 
respectively (it is assumed that $1\le r< n$), 
 such that every    point $p\in M^n$
has  a   neighborhood 
$U´(p)$   with   
coordinates  
$$
(\bar x,\bar y)= \bigr((x_1,x_2,...x_r),(y_{r+1},y_{r+2},...,y_n)\bigl)
$$
such that  the  $x$-coordinates are 
constant on every leaf  of the foliation $B_{n-r}\cap U´(p)$, 
the $y$-coordinates are 
constant on every leaf of the foliation $B_{r}\cap U´(p)$,   and 
 the metric  $h$ is block-diagonal   such that the first 
($r\times r$) block depends on the $x$-coordinates and 
the last  $((n-r)\times (n-r))$ block depends on the $y$-coordinates. 
\end{Def}

A model example of manifolds with local-product structure is 
the direct product of two Riemannian manifolds $(M_1^r, g_1)$ and   
$(M_2^{n-r}, g_2)$.  In this case, the leaves  of the foliation 
$B_r$ are the products  of $M_1^r$ and  the points of $M_2^{n-r}$, 
the leaves  of the foliation 
$B_{n-r}$ are the products  of  the points of $M_1^{r}$ and  $M_2^{n-r}$, and
the  metric $h$ is the product metric $g_1+g_2$.

Below we assume that 
\begin{itemize}
\item[(a)] The metrics $g$ and $\bar g$ are projectively equivalent 
on a connected  $M^n$. 
\item[(b)]
 They 
are strictly-non-proportional at least
 at one point of $M^n$.
                     \item[(c)]   There exists $r$,  
$1\le r<n$,  such  that  
$\lambda_r<\lambda_{r+1}$ at every  point of $M^n$. 
\end{itemize} 

We will show that (under the assumptions (a,b,c)) we can naturally define 
two local-product structures $(h, B_r, B_{n-r})$ and $(\bar h, B_r, B_{n-r})$
such that  the restrictions of $h$  and
$\bar h$ to every  leaf are  projectively equivalent and 
strictly non-proportional at least at one point.

 At every 
 point $x\in M^n$, denote by $V^r_x$ 
  the  subspaces  of $T_xM^n$ spanned by the eigenvectors of 
$L$ corresponding to the eigenvalues 
$\lambda_1,...,\lambda_r$. Similarly,  denote by $V^{n-r}_x$ 
  the  subspaces  of $T_xM^n$ spanned by the eigenvectors of 
$L$ corresponding to the eigenvalues 
$\lambda_{r+1},...,\lambda_n$. By assumption, for any $i,j$ such that 
 $i\le r < j$, we have $\lambda_i\ne \lambda_j$ so that 
$V^r_x$ and $V^{n-r}_x$ are  two  
 smooth distributions  on $M^n$. By Corollary~\ref{nijenhuis}, 
the distributions are integrable so that they define two 
transversal     foliations $B_r$ and $B_{n-r}$ of dimensions $r$ and $n-r$, 
respectively.

By construction,   the distributions $V_r$ and 
 $V_{n-r}$ are invariant with respect to $L$.  Let us denote by 
$L_r$,  $L_{n-r}$ the restrictions of $L$ to $V_r$ and 
 $V_{n-r}$, respectively.  
We will  
denote by $\chi_r$, $\chi_{n-r}$  the characteristic polynomials of
$L_r$,  $L_{n-r}$, respectively. 

 Consider the 
 (1,1)-tensor 
$$   
C\eqdef  \left((-1)^{r}\chi_r(L)+\chi_{n-r}(L)\right)
$$ 
and the metric $h$ given by the relation
$$  
h(u,v)\eqdef g(C^{-1}(u),v)
$$ 
for any  vectors $u,v$. (In the 
 tensor notations, the metrics  $h$ and $g$ are related by 
$g_{ij}=h_{i\alpha}C_j^\alpha$). 

 Consider the 
 (1,1)-tensor 
$$   
\bar C\eqdef  \left(\frac{(-1)^{r}}{\det(L_{n-r})}\chi_r(L)+\frac{1}{\det(L_{r})}\chi_{n-r}(L)\right)
$$ 
and the metric $\bar h$ given by the relation
$$  
\bar h(u,v)\eqdef \bar g(\bar C^{-1}(u),v)
$$ 
for any  vectors $u,v$.

\begin{Lemma}\hspace{-2mm}{\bf .}
 \label{local-direct-product} The following statements hold: 

 \begin{enumerate}
  \item The triples $(h, B_r, B_{n-r})$ and $(\bar h, B_r, B_{n-r})$ 
  are local-product structures on $M^n$.

  \item   For any leaf of $B_r$, the restrictions of $h$  and
$\bar h$ to the leaf are  projectively equivalent and 
strictly non-proportional at least at one point. For any leaf of $B_{n-r}$,
 the restrictions of $h$  and
$\bar h$ to the leaf are  projectively equivalent and 
strictly non-proportional at least at one point. 
\end{enumerate}
\end{Lemma}

\noindent {\bf Proof:}  First of all,  $h$ and $\bar h$  are well-defined
 Riemannian metrics.
Indeed, take an arbitrary point $x\in M^n$. At the tangent space to this point, 
  we can find a coordinate system where the tensor $L$  and the metric $g$ 
are diagonal.  In this coordinate system, 
the characteristic polynomials  $\chi_r$, $\chi_{n-r}$ are given by 
\begin{equation}\label{lcr}
\begin{array}{lcr}
(-1)^{r}\chi_r&=& (t-\lambda_1)(t-\lambda_2)...(t-\lambda_r)   \\
\chi_{n-r}&=& (\lambda_{r+1}-t)(\lambda_{r+2}-t)...(\lambda_n-t). 
\end{array}
\end{equation}
Then,  the (1,1)-tensors 
\begin{eqnarray*}
C&=&\left((-1)^{r}\chi_r(L)+\chi_{n-r}(L)\right)  \\
\bar C&=&\left({(-1)^{r}}{\det(L_{n-r})}\chi_r(L)+{\det(L_{r})}\chi_{n-r}(L)\right)
\end{eqnarray*}
\noindent  are  given by 
the  diagonal matrices  
\begin{eqnarray} 
\hspace{-5mm}\diag\left(\prod_{j=r+1}^n(\lambda_{j}-\lambda_1),..., 
      \prod_{j=r+1}^n(\lambda_{j}-\lambda_r),                              
 \prod_{j=1}^r(\lambda_{r+1}-\lambda_j),...,
      \prod_{j=1}^r(\lambda_{n}-\lambda_j)\right), \label{chi1} \\
\nonumber
 \hspace{-5mm}\diag\left(\frac{1}{\lambda_{r+1}...\lambda_{n}}\prod_{j=r+1}^n(\lambda_{j}-\lambda_1),..., 
     \frac{1}{\lambda_{r+1}...\lambda_{n}} \prod_{j=r+1}^n(\lambda_{j}-\lambda_r),\right.\\  \left.
 \frac{1}{\lambda_{1}...\lambda_{r}}                             
 \prod_{j=1}^r(\lambda_{r+1}-\lambda_j),...,
     \frac{1}{\lambda_{1}...\lambda_{r}} \prod_{j=1}^r(\lambda_{n}-\lambda_j)\right). \label{chi2}
\end{eqnarray} 
We see that the tensors are   diagonal and 
 all its  diagonal 
components
  are positive.  
Then,  the tensors  $C^{-1}, \bar C^{-1}$ are  well-defined and 
$h,\bar h$ are   Riemannian metrics.

By construction, $B_r$ and $B_{n-r}$ are well-defined 
transversal foliations of supplementary 
dimensions.  In order to prove Lemma~\ref{local-direct-product}, 
we need to verify 
 that,  locally,  the triples  $(h,B_r,B_{n-r})$
 and $(\bar h,B_r,B_{n-r})$  are as in Definition~\ref{lp},  that the 
restriction of $h$   to a leaf is  projectively equivalent to the restriction of $\bar h$, 
and that  the 
restriction of $h$   to a leaf  is strictly non-proportional 
(at least at a point)
to the restriction of $\bar h$.    

It is  sufficient to verify the first two  statements 
 at almost every point of $M^n$. More precisely,
it is known that 
the 
 triple $(h, B_r, B_{n-r})$ is a local-product 
structure if and only if the foliations $B_r$ and 
$B_{n-r}$ are orthogonal and totally geodesic \cite{Naveira}.
 Clearly, if the foliations and the metric are globally given and smooth, 
 if the foliations are orthogonal and totally-geodesic at 
almost every point then they   are orthogonal and totally-geodesic 
at every point. 

Similarly,  if a   foliation and two  metrics are globally-given and smooth, 
if the restriction of the metrics to the leaves of the foliation 
is projectively equivalent  almost everywhere
 then it is so at every point.

By Corollary~\ref{ordered12}, at 
almost every point of $M^n$ the eigenvalues of $L$ are different. 
Consider the 
 Levi-Civita coordinates  
  $x_1,..., x_n$ from Theorem~\ref{Levi-Civita}.
 In the Levi-Civita  coordinates, $L$ is given by (\ref{lambda}). 
Then, by constructions of the foliations $B_r$ and $B_{n-r}$, 
 the coordinates 
$x_1,...,x_r$
  are  constant on every leaf  of the foliation $B_{n-r}$,  
the coordinates $x_{r+1},...,x_n$ 
 are 
constant on every leaf of the foliation $B_{r}$. 

Using (\ref{chi1},\ref{chi2} ), we see that,  
in the Levi-Civita coordinates, 
$h, \bar h$ are  given by 
\begin{eqnarray}
h(\dot{x}, \dot{x})&=&\ \tilde \Pi_1dx_1^2+...+
\ \tilde \Pi_rdx_r^2+\nonumber\\
&+&\tilde \Pi_{r+1}dx_{r+1}^2+...+\tilde \Pi_ndx_n^2, \label{h}\\
\bar h(\dot{x}, \dot{x})&=&\ {\tilde\rho_1}\tilde \Pi_1dx_1^2+\cdots+
\ {\tilde\rho_r}\tilde\Pi_rdx_r^2+\nonumber\\
&+&{\tilde\rho_{r+1}}\tilde\Pi_{r+1}dx_{r+1}^2 \cdots {\tilde\rho_n}
\tilde\Pi_ndx_n^2, \label{bh}
\end{eqnarray}
where the functions $\tilde \Pi_i,\tilde \rho_i$ 
are as follows: 
for   $i\le r$, they  are  given by 
\begin{eqnarray*}
\tilde \Pi_i &\eqdef &
{(\lambda_i-\lambda_1)...(\lambda_i-\lambda_{i-1})(\lambda_{i+1}-\lambda_i)...(\lambda_r-\lambda_i)}, \\
\tilde\rho_i &\eqdef& \frac{1}{\lambda_i (\lambda_1\lambda_2...\lambda_r)}.
\end{eqnarray*}
For   $i> r$,   the functions $\tilde \Pi_i,\tilde \rho_i$ are  given by 
\begin{eqnarray*}
\tilde \Pi_i & \eqdef & 
{(\lambda_i-\lambda_{r+1})...(\lambda_i-\lambda_{i-1})(\lambda_{i+1}-\lambda_i)...(\lambda_n-\lambda_i)}, \\
\tilde\rho_i  & \eqdef &  \frac{1}{\lambda_i (\lambda_{r+1}\lambda_2...\lambda_{n})}. 
 \end{eqnarray*}

We see that the restrictions of the metrics on the leaves of the 
foliations have the form from Levi-Civita's Theorem and, therefore, 
are projectively equivalent. We see that the metrics $h$, $\bar h$ are block-diagonal
with the first $r\times r$ block depending on the variables $x_1,...,x_r$ and 
the second $(n-r)\times (n-r)$ block depending on the remaining variables, so that 
$(h, B_r, B_{n-r})$ and $(\bar h, B_r, B_{n-r})$ are local-product structure. 

The last thing to show is that the restrictions of the metrics to  every  leaf
are strictly non-proportional at least at one point. Suppose it is not 
so; that is, there exists a leaf (say, of foliation $B_r$) and $k$, $0<k<r$ 
such that $\lambda_k=\lambda_{k+1}$ at each point of the leaf. Then, by 
the first statement of Corollary~\ref{ordered12}, the eigenvalues 
 $\lambda_k$, $\lambda_{k+1}$ are constant on the leaf. Since the Nijenhuis torsion of $L$ is zero, the eigenvalues $\lambda_k$, $\lambda_{k+1}$ are 
constant 
along the leaves of the foliation $B_{n-r}$ \cite{Haantjes}. Then 
$\lambda_k=\lambda_{k+1}$ at each point of a neighborhood of the leaf, which 
contradicts  Corollary~\ref{multiplicity}.  
Thus the restrictions of the metrics to  every  leaf
are strictly non-proportional at least at one point. 
  Lemma~\ref{local-direct-product}
is proved.

\subsection{Gluing procedure} \label{2.3}

Let $B_r$ and $B_{n-r}$, $0<r<n$, be transversal foliations of 
supplementary dimensions $r$ and $n-r$ on a connected  
$M^n$.  Suppose there exist Riemannian metrics 
$h$ and $\bar h$ such that 
\begin{itemize} 
\item[(i)] The triples $(h, B_r, B_{n-r})$ and $(\bar h, B_r, B_{n-r})$
are local-product structures. 

\item[(ii)] For every fiber of $B_r$ the restrictions of $ h$, $\bar h$
to the fiber are projectively equivalent and strictly-non-proportional 
at least at one point. For every fiber of $B_{n-r}$ 
the restrictions of $ h$, $\bar h$
to the fiber are projectively equivalent and strictly-non-proportional 
at least at one point.
\end{itemize} 

Let us denote by $L_r$ (by $L_{n-r}$, respectively) 
 the tensor (\ref{l}) constructed for the restrictions of   
$h$, $\bar h$ to the tangent spaces of the leaves of  
$B_r$ (of $B_{n-r}$, respectively). Assume in addition that 
\begin{itemize}
\item[(iii)] for any point of $M^n$, the eigenvalues of 
 $L_{n-r}$ are greater than that of 
$L_{r}$. 
\end{itemize}

At every point of $M^n$, let us denote by 
 $\chi_r$, $\chi_{n-r}$  the characteristic polynomials of
$L_r$,  $L_{n-r}$, respectively. 
Let us denote by $P_r:T_xM^n\to T_xB_r$,
$P_{n-r}:T_xM^n\to T_xB_{n-r}$   
 the orthogonal  projections to the tangent spaces of the 
foliations $B_r$, $B_{n-r}$. (Since the foliations are 
orthogonal with respect to both metrics, it does not matter 
what metric we take here.) There exists a unique $(1,1)$-tensor $L$
 such that 
 $L\circ P_r = L_r$ and  $L\circ P_{n-r} = L_{n-r}$.

 Consider the 
 (1,1)-tensor 
\begin{equation}\label{C}   
C\eqdef  \left((-1)^{r}\chi_r(L)+\chi_{n-r}(L)\right)
\end{equation}
and the metric $g$ given by the relation
\begin{equation}\label{g}    
g(u,v)\eqdef h(C(u),v)
\end{equation}
for any  vectors $u,v$. (In the 
 tensor notations, the metric $g$ is given by
$h_{i\alpha}C_j^\alpha$.)

 Consider the 
 (1,1)-tensor 
   \begin{equation}\label{bC}  
\bar C\eqdef  \left(\frac{(-1)^{r}}{\mbox{det}(L_{n-r})}\chi_r(L)+\frac{1}{\mbox{det}(L_{r})}\chi_{n-r}(L)\right)
\end{equation}
and the metric $\bar h$ given by the relation
\begin{equation}\label{bg}
\bar g(u,v)\eqdef \bar h(\bar C(u),v)
\end{equation}
for any  vectors $u,v$.

\begin{Lemma}\hspace{-2mm}{\bf .} \label{inverse} 
The metrics $g$ and $\bar g$ are projectively equivalent on 
$M^n$  and are strictly-non-proportional at almost every 
 point of $M^n$. 
\end{Lemma}

\noindent {\bf Proof:}  First of all,  $g$ and $\bar g$ 
are well-defined Riemannian metrics. More 
precisely, 
at the tangent space to every point we can 
find a coordinate system where $L$ is given by the 
diagonal matrix $\diag(\lambda_1,...,\lambda_n)$ assuming 
$$
\lambda_1\le \lambda_2\le ...\le\lambda_r<\lambda_{r+1}\le ...\le\lambda_n.
$$
In this coordinate system, 
 the characteristic polynomials $\chi_r$ and $\chi_{n-r}$
are given by (\ref{lcr}). Hence, the tensors 
$C,\bar C$ are  given by (\ref{chi1},\ref{chi2}), and the condition (iii) 
guarantees that 
the metrics $g$ and $\bar g$ are well-defined Riemannian metrics. 

Let us show that the metrics  $g$ and $\bar g$ are projectively equivalent. 
It is sufficient to check it at almost every point. Clearly, at 
almost every 
point  of $M^n$ the eigenvalues of the tensor $L$ are all different. 
Since by  Levi-Civita's Theorem, the restriction of 
$h$, $\bar h$ to the leaves  of $B_r$, $B_{n-r}$ 
have the model form (\ref{Canon_g},\ref{Canon_bg}), 
we  obtain that the metric $h, \bar h$ are given by (\ref{h},\ref{bh}). 
Using that $C,\bar C$ are  given by (\ref{chi1},\ref{chi2}), we obtain that 
the metrics $g$, $\bar g$ have precisely the form from Levi-Civita's Theorem,
 and, therefore, are projectively equivalent. Lemma~\ref{inverse} is proved.

\begin{Rem}\hspace{-2mm}{\bf .} \label{missleading} 
The notation $C$, $\bar C$, $L$,   $\lambda_i$, $g$ and $\bar g$ used in this section  are not misleading and are 
compatible with the notations in Section~\ref{2.2}. More precisely, 
 if we take $g, \bar g$ satisfying assumptions 
(a,b,c), construct the metrics $h, \bar h$ and foliations $B_r, B_{n-r}$, 
then  (by Lemma~\ref{local-direct-product}) the triples 
$(h, B_r, B_{n-r})$,  $(\bar h, B_r, B_{n-r})$ satisfy conditions 
(i, ii, iii).  Moreover,  the tensor $L$ given by (\ref{l}) coincide with 
the tensor $L$   constructed in this section. 
Therefore, 
the tensors (\ref{C},\ref{bC}) coincide with the 
tensors $C,\bar C$ from Section~\ref{2.2}, and, therefore, 
the metrics constructed by (\ref{g},\ref{bg}) coincide 
with the initial metrics $g,\bar g$. 
\end{Rem}

\begin{Rem}\hspace{-2mm}{\bf .} \label{example}
Levi-Civita's Theorem~\ref{Levi-Civita}
follows from  Lemmas~\ref{local-direct-product},\ref{inverse}. 
\end{Rem}
\noindent {\bf Proof of Remark \ref{example}:} 
First of all, by direct calculation it is 
possible to verify that if the metrics are given by the Levi-Civita's model 
form  (\ref{Canon_g},\ref{Canon_bg}), then they are projectively equivalent. 

In order  to prove that strictly non-proportional projectively equivalent 
 metrics have (locally)  the form  (\ref{Canon_g},\ref{Canon_bg}), we use induction in dimension of the manifold. 

 For one-dimensional manifold it is 
nothing to prove; suppose Levi-Civita's Theorem is true 
for dimension $n-1$. 
Let us prove that, for dimension $n$, strictly non-proportional 
projectively equivalent metrics are locally given by (\ref{Canon_g},\ref{Canon_bg}). 
 
If  the metrics are 
 strictly non-proportional at $p$, then $\lambda_1<\lambda_2$ in a small neighborhood of $p$. 
Put $r=1$ and construct the local-direct-product structures
$(h, B_1, B_{n-1})$ and $(\bar h, B_1, B_{n-1})$.
By definition~\ref{lp}, 
 there exists a smaller neighborhood  of $p$ such that the foliations $B_1$ and  $B_{n-1}$ there look as they came from the direct product of the interval 
and the (n-1)-dimensional disk. Let us choose on leaf  of the 
foliation $B_1$ and one leaf  of the foliation $B_{n-1}$.  
 
Since the leaf of $B_1$ is one-dimensional, there exists a coordinate $x_1$
 there such that the restriction of the metrics $h, \bar h$ are respectively 
given by $dx_1^2$ and $\frac{dx_1^2}{\lambda_1^2}$. Since the restrictions of
 the 
metrics $h, \bar h$ to the leaf of $B_{n-1}$ are strictly non-proportional and 
projectively equivalent, by Levi-Civita's Theorem, there exists 
a coordinate system $x_2,...,x_n$ there such that the  restrictions of the 
metrics $h, \bar h$ to the leaf of $B_{n-1}$ are respectively given by 
\begin{eqnarray*}
\hat \Pi_2 dx^2_2&+...+&  \hat \Pi_n dx^2_n\\
\frac{1}{\hat \rho_2}\hat \Pi_2 dx^2_2&+...+& \frac{1}{\hat \rho_n} \hat \Pi_n dx^2_n, 
\end{eqnarray*}
where the functions $\hat \rho_i$ and $\hat \Pi_i$  are related to the 
functions $\rho_i$ and $\Pi_i$  from 
Levi-Civita's Theorem by the formulae
$$
\hat\rho_i=\lambda_1\rho_i, \ \ \ \hat\Pi_i=\frac{\Pi_i}{\lambda_i-\lambda_1}.
$$
Because of the local-product structure, these coordinates of the leaf of $B_1$ and on the leaf of  $B_{n-1}$ give us a coordinate system in the neighborhood of $p$. 
 By direct calculations, $-\chi_1=(t-\lambda_1)$, $\chi_{n-1}=
(\lambda_2-t)(\lambda_3-t)...(\lambda_n-t)$. Then, 
$$
-\chi_1(L)=\diag(0,\lambda_2-\lambda_1, \lambda_3-\lambda_1,...,\lambda_n-\lambda_1),
$$
$$
\chi_{n-1}(L)=\diag(\Pi_1,0,0,...,0), 
$$
and the determinants $\mbox{det}(L_{1}),\mbox{det}(L_{n-1})$ are equal to  $\lambda_1$, 
$\lambda_2\lambda_3...\lambda_n$, respectively.  

Using that the metric $h$, $\bar h$ are the products of 
their  restrictions to the leaf of $B_1$ and the leaf of $B_{n-1}$, 
and in view of Remark~\ref{missleading}, 
 we obtain that the metrics $g, \bar g$  are precisely in the model form 
(\ref{Canon_g},\ref{Canon_bg}).

\subsection{Dimension 2} \label{2.4} 

\noindent The goal of this section is to give the local description
             of  projectively equivalent metrics (on surfaces) 
             near the points where the eigenvalues of $L$ bifurcate. 
           In dimension two,  the inverse of Theorem~\ref{integrability} also 
          takes place:   
\begin{Th}[\cite{MT,quantum,surface}]\label{dim2}
Let $g$, $\bar g$ be Riemannian metrics on $M^2$. 
They are projectively equivalent if and only if the 
function 
\begin{equation} \label{integral2}
F:TM^2\to R,
\ \ F(\xi)\eqdef 
\left(\frac{\mbox{det}(g)}{\mbox{det}(\bar g)}\right)^{\frac{2}{3}} 
 \bar g(\xi,\xi)
\end{equation}
is an integral of the geodesic flow of $g$. 
\end{Th}

\noindent We see that the integral $F$ is quadratic in velocities. 
          Thus the existence of an integral quadratic in 
          velocities (for the geodesic flow of $g$) 
          allows one to construct a metric projectively equivalent to $g$
  (at least locally). 
         
         Now, in   
         in the two-dimensional case, the local description of metrics with 
         quadratically integrable geodesic flows has been obtained  
in \cite{BMF}, see also \cite{BF}  (Basing on the technique developed in \cite{Kol}). Combining this description with Theorem~\ref{dim2}, we obtain the 
following

 \begin{Th}\label{local2}\hspace{-2mm}{\bf .}    
  Let $g$ and $\bar g$ be projectively equivalent on a 
(2-dimensional)  connected surface $M^2$. Suppose they are 
non-proportional at least at one point. Assume they  are proportional at 
$p\in M^2$.  Then, precisely one of the following possibilities takes place:

\begin{enumerate}
\item   
 There exist coordinates $u,v$ in a neighborhood of $p$,  and
there exists  function $\lambda$ of one variable 
such that the metrics have the following model form  
 \begin{eqnarray} 
 ds_g^2 & =  & 2\frac{\lambda(\rho+u)-\lambda(u-\rho)}{\rho}(du^2+dv^2) \label{mg}\\
  ds_{\bar g}^2 & =  & \nonumber
\left(\frac{\lambda(u+\rho)-\lambda(u-\rho)}{\rho\lambda(u-\rho)\lambda(u+\rho)}\right)^2
\left[\left(\rho\frac{\lambda(u+\rho)+\lambda(u-\rho)}{\lambda(u+\rho)-\lambda(u-\rho)} \right)(du^2+dv^2)\right.\\
&-& \left. udu^2-2vdudv+udv^2\right]
\label{mbg}
  \end{eqnarray}
where 
$\rho\eqdef \sqrt{u^2+v^2}$. 
 
\item
There exist  coordinates $u,v$ in a neighborhood of $p$,  
there exists   a  functions $f$ of one variable, and there exists a positive 
constant $\lambda_1$  
such that the metrics have the following model form

\begin{eqnarray} 
 ds_g^2 & =  & {f(\rho^2)}(du^2+dv^2)\label{l0} \\ \label{l1}
ds_{\bar g}^2 & = & 
\frac{f(\rho^2)}{\lambda_1(\lambda_1+\lambda_1\rho^2f(\rho^2))^2}\left(\bigl(1+f(\rho^2)v^2\bigr)
du^2- 2f(\rho^2)uvdudv\right. \nonumber \\  
& + &\left.\bigl(1+f(\rho^2)u^2\bigr)dv^2\right),
\end{eqnarray}
where 
$\rho\eqdef \sqrt{u^2+v^2}$.

\item
There exist  coordinates $u,v$ in a neighborhood of $p$,  
there exists    function $f$ of one variable, and there exists a positive 
constant $\lambda_2$  
such that the metric  $g$  has the form (\ref{l0})
and the metric $\bar g$  is given by

\begin{eqnarray}  \label{l2}
ds_{\bar g}^2 & = & \frac{f(\rho^2)}{\lambda_2(\lambda_2-\lambda_2f(\rho^2))^2}\left(\bigl(1-f(\rho^2)v^2\bigr)du^2+ 2f(\rho^2)uvdudv\right. \nonumber \\  & + &\left.\bigl(1-f(\rho^2)u^2\bigr)dv^2\right),
\end{eqnarray}
where 
$\rho\eqdef \sqrt{u^2+v^2}$.
\end{enumerate}  
 \end{Th} 

Theorem~\ref{local2} is true also in the other direction: 
if Riemannian  metrics are 
given by formulae   (\ref{mg},\ref{mbg}) or  (\ref{l0},\ref{l1}) or 
 (\ref{l0},\ref{l2}), then          they are projectively equivalent. 
Of course, in order the formulae to  define 
  Riemannian  metrics (at least in a small neighborhood of $(0,0)$), 
the functions $f$ and $\lambda$ must be smooth, 
positive and satisfy the conditions $\lambda'(0)>0$, $f'(0)>0$.

\begin{Rem} \hspace{-2mm}{\bf .} The most natural coordinate system 
for projectively equivalent metrics near a point of bifurcation has 
singularity at the point: 
\begin{enumerate} 
\item In the  elliptic coordinate system  
$x^2_1=\rho-u$, $x^2_2=\rho+u$
(which  has  a singularity at $(0,0)$), the metrics 
(\ref{mg},\ref{mbg})
are given (up to multiplication by 4) by  
$$
(\lambda(x_1)-\lambda(x_2))(dx_1^2+dx_2^2),
$$
$$
\left(\frac{1}{\lambda(x_2)}-\frac{1}{\lambda(x_1)}\right)
\left(\frac{dx_1^2}{\lambda(x_1)}+\frac{dx_2^2}{\lambda(x_2)}\right),  
$$ which is precisely the Levi-Civita form for dimension two. 
In particular, the eigenvalues $\lambda_1(u,v)$,  $\lambda_2(u,v)$ of $L$  are equal to  $\lambda(u-\rho)$, $\lambda(u+\rho)$, respectively. 

\item 
In the polar coordinate 
system $u=e^{r} \cos(\phi), v=e^r\sin(\phi)$ the metrics 
 (\ref{l0},\ref{l1}) and 
(\ref{l2}) are given by 
 \begin{eqnarray*}
 (e^{2r}f(e^{2r}))\bigl(dr^2 &+& d\phi^2\bigr), 
\\
\frac{1}{\lambda_1}\left(\frac{1}{\lambda_1}-\frac{1}{\lambda_1+\lambda_1e^{2r}f(e^{2r})}\right)  
\left( \frac{d\phi^2}{\lambda_1} \right.& + &\left. \frac{dr^2}{\lambda_1+\lambda_1e^{2r}f(e^{2r})}
\right),
 \\
\frac{1}{\lambda_2}\left(\frac{1}{\lambda_2-\lambda_2e^{2r}f(e^{2r})}-\frac{1}{\lambda_2}\right)
\left(\frac{dr^2}{\lambda_2-\lambda_2e^{2r}f(e^{2r})}\right. & + & \left.
\frac{d\phi^2}{\lambda_2}\right),  
\end{eqnarray*}
respectively. 
We see that they are in the Levi-Civita form (up to the factors $\lambda_1$ and $\lambda_2$), and that
the eigenvalues of $L$ for the pair of metrics  
(\ref{l0},\ref{l1}) are  
$\lambda_1(u,v) = \lambda_1$, $\lambda_2(u,v) = 
\lambda_1+\lambda_1\rho^2f(\rho^2)$, and 
the eigenvalues of $L$ for the pair of metrics  (\ref{l0},\ref{l2}) are 
$\lambda_1(u,v) = \lambda_2-\lambda_2\rho^2f(\rho^2)$,
$\lambda_2(u,v) = \lambda_2$.
\end{enumerate}
\end{Rem}
  
Wee see that the first possibility 
(from Theorem~\ref{local2}) 
 for projectively equivalent metrics 
is characterized by the  condition that
$\lambda_1$,  $\lambda_2$ of  $L$  for 
  are non-constant, and the second   possibility   is characterized by the 
condition $\lambda_1$ is constant, $\lambda_2$ is not constant, and the third 
   possibility   is characterized by the 
condition $\lambda_2$ is constant, $\lambda_1$ is not constant.

\subsection{Dimension 3} \label{2.5}
The goal of this section is to describe 
(strictly-non-proportional at least at a  point)
projectively equivalent metrics on a 3-manifold near the points where 
the metrics are  proportional. This will be made in Theorem~\ref{local3}. 
In order to prove it, we need Corollary~\ref{l3} and Lemma~\ref{fixed}.

\begin{Cor} \hspace{-2mm}{\bf .} \label{l3}
Let $g$ and $\bar g$ be projectively equivalent on connected 
$M^3$ and be strictly non-proportional at least at one point of 
$M^3$.  Let $p\in M^3$. Suppose  in a neighborhood 
of $p\in M^3$ the eigenvalue  $\lambda_2$ is constant, and suppose 
$\lambda_1(p)<\lambda_2=\lambda_3(p)$.
Then, there exists a neighborhood of $p$ with coordinates $x_1,x_2,x_3$ where
the metrics have the following model form: 
\begin{eqnarray}
ds_g^2&=& (\lambda_1(x_1)-1)(\lambda_1(x_1)-1-\rho^2f(\rho^2))dx_1^2+
f(\rho^2)\left(1+x_2^2f(\rho^2)-\lambda_1(x_1)\right) dx_2^2   \label{l4}
\nonumber \\
&-& 2f(\rho^2)x_2x_3dx_2^2dx_3^2 + \left(1+x_3^2f(\rho^2)-\lambda_1(x_1)\right)dx_3^2,\\  
ds_{\bar g}^2&=& 
\frac{(\lambda_1(x_1)-1)(\lambda_1(x_1)-1-\rho^2f(\rho^2))}{\lambda_1^2(x_1)(1+
\rho^2f(\rho^2))}dx_1^2\nonumber\\&+&\frac{f(\rho^2)}{\lambda_1(x_1)(1+\rho^2f(\rho^2))}
\left[\left(1-\lambda_1(x_1)\frac{1+x_3^2f(\rho^2)}{1+\rho^2f(\rho^2)}\right)dx_2^2\right.\nonumber\\
&+&
\frac{2\lambda_1(x_1)x_2x_3}{1+\rho^2f(\rho^2)}dx_2dx_3\nonumber\\
&+&\left.\left(1-\lambda_1(x_1)\frac{1+x_2^2f(\rho^2)}{1+\rho^2f(\rho^2)}\right)dx_3^2\right], \label{l5}
\end{eqnarray} 
 where $\rho\eqdef \sqrt{x_2^2+x_3^2}$; 
 $f$ and $\lambda_1$ are 
functions of one  variable and $\lambda_2$ is a positive constant.
\end{Cor}

\begin{Rem}\hspace{-2mm}{\bf .} 
In the cylindrical  coordinates  $x_1=u_1,x_2=u_3 \cos(u_2), v=u_3 \sin(u_2)$
 the 
metrics $g$, $\bar g$ almost have Levi-Civita  form 
(\ref{Canon_g},\ref{Canon_bg}).  
\end{Rem} 

Proof of Corollary~\ref{l3}: Since $\lambda_1<\lambda_2$ at the point $p$, 
there exists a neighborhood of $p$ where  $\lambda_1<\lambda_2$. Put $r=1$,
apply the splitting procedure from Section~\ref{2.2} and construct the
metrics $h$ and $\bar h$ and the  
foliations  $B_1$, $B_2$. 
 By 
Lemma~\ref{local-direct-product}, there exists 
 a (possible, smaller) neighborhood of $p$ isomorphic to a direct 
 product of an interval and a 2-disc such that the 
metric $h$ is the product metric $g_1+g_2$ (where the metric $g_1$
is a metric on the interval and $g_2$ is a metric on the disc) 
and the metric $\bar h$ is  the product metric $\bar g_1+\bar g_2$
(where the metric $\bar g_1$
is a metric on the interval and $\bar g_2$  is a metric on the disc
projectively equivalent to the metric $g_2$ and strictly-non-proportional to $g_2$
at least at one point). 
Since $\lambda_2$ is constant,  the smallest eigenvalue of the tensor $L$ 
constructed for the metrics $g_2$, $\bar g_2$  is constant. Since 
$\lambda_2=\lambda_3$ at $p$,  
 the 
metrics $g_2$, $\bar g_2$ are proportional at one point and, 
 therefore, 
 are  given by (\ref{l0},\ref{l1}) in an appropriate coordinate system $x_2, 
x_3$.
There evidently exists a coordinate $x_1$ on the interval such that the 
metrics $g_1$, $\bar g_1$ are given by $dx_1^2$, $\frac{dx_1^2}{\lambda_1^2(x_1)}$
for  an appropriate function $\lambda$.  Applying the gluing procedure from 
Section~\ref{2.3}, we obtain precisely the form (\ref{l4},\ref{l5}). 
Corollary~\ref{l3} is proved. 

\begin{Lemma}\hspace{-2mm}{\bf .} Consider the Riemannian metrics \label{fixed} 
given by the 
formulae (\ref{l4},\ref{l5}) in a neighborhood of the point $(0,0,0)$. 
 Then, 
\begin{itemize} 
\item  the plane $P:=\{(x_1,x_2,x_3):\ x_2=0\}$ 
is a totally geodesic submanifold; 
\item the eigenvalue $\lambda_3$  equals $\lambda_2$ precisely at the points where 
$x_2=x_3=0$; 
\item  the action of $S^1$ by the rotations 
 $$(x_1,x_2,x_3)\mapsto (x_1, x_2\cos(\phi)-x_3\sin(\phi), 
x_2\sin(\phi)+x_3\cos(\phi))$$ preserves both metrics,
 and at every   point  $x\not\in \{(x_1,x_2,x_3):\ x_2=x_3=0\}$
 its orbits are tangent to the eigenspace of $L$ corresponding to $\lambda_2$, 
\item  at every point of the plane $P$, 
the vector $\left(\frac{\partial }{\partial x_2}\right)$ is the eigenvector
of $L$ with the eigenvalue $\lambda_2$.  
\end{itemize}
\end{Lemma}

The lemma could  be proved by direct computations. Actually, the first
 statement follows from the fact that the symmetry 
$(x_1,x_2,x_3)\mapsto (x_1,-x_2,x_3)$ is evidently an isometry; the second,
 third  and the fourth  statements follow
 from the observation that splitting-gluing
 procedure using in 
construction of metrics (\ref{l4},\ref{l5}) is invariantly given in terms of metrics and therefore inherits all symmetries of the metrics $g_1$ and $g_2$.

\begin{Th}\hspace{-2mm}{\bf .}\label{local3}    
  Let $g$ and $\bar g$ be projectively equivalent on a 
(3-dimensional)  connected manifold  $M^3$. Suppose they are strictly
non-proportional at least at one point. Assume they  are proportional at 
$p\in M^3$.  Then, 
 there exist coordinates $u_1,u_2,u_3$ in a neighborhood of $p$,  
a   function $\lambda$ of one variable  and a positive constant $C$,  
such that the metrics have the following model form  
 \begin{eqnarray} 
 ds_g^2 & =  & 2\frac{\lambda(\rho+u_1)-\lambda(u_1-\rho)}{\rho}\left(du_1^2+\left(d\sqrt{u_2^2+u_3^2}\right)^2\right)\nonumber \\ &+&
C(\lambda(0)-\lambda(u_1-\rho))(\lambda(u_1+\rho)-\lambda(0))\left(u_3du_2-u_2du_3\right)^2 
\label{mg3}\\
  ds_{\bar g}^2 & =  & \nonumber
\left(\frac{\lambda(u_1+\rho)-\lambda(u_1-\rho)}{\rho\lambda(u_1-\rho)\lambda(u_1+\rho)}\right)^2
\left[\left(\rho\frac{\lambda(u_1+\rho)+\lambda(u_1-\rho)}{\lambda(u_1+\rho)-\lambda(u_1-\rho)} \right)\left(du_1^2+\left(d\sqrt{u_2^2+u_3^2}\right)^2\right)\right.\\
&-& \left. u_1du_1^2-2\sqrt{u_2^2+u_3^2}du_1d\sqrt{u_2^2+u_3^2} +u_1\left(d\sqrt{u_2^2+u_3^2}\right)^2\right]\nonumber \\
&+&C\frac{(\lambda(0)-\lambda(u_1-\rho))(\lambda(u_1+\rho)-\lambda(0))}{\lambda(0)^2\lambda(u_1-\rho)\lambda(u_1+\rho)}\left(u_3du_2-u_2du_3\right)^2 
\label{mbg3}
  \end{eqnarray}
where 
$\rho\eqdef \sqrt{u_1^2+u_2^2+u_3^2}$. 
 \end{Th}

\begin{Rem}\hspace{-2mm}{\bf .} The most natural coordinate system
here are  cylindrical-elliptic: 
$$
x_1=\rho-u_1; \ x_2=\sqrt{C}\arccos\left(\frac{u_2}{\sqrt{u_2^2+u_3^2}}\right); \ x_3=\rho+u_1, 
$$
where the metrics have the Levi-Civita form (\ref{Canon_g},\ref{Canon_bg})
(with $\lambda_1=\lambda(\rho-u_1)$, $\lambda_2=\lambda(0)$, and $\lambda_3= \lambda(\rho+u_1)$.) In particular, if two Riemannian 
metrics are given in the form (\ref{mg3},\ref{mbg3}), they are projectively equivalent. If $\lambda$ is smooth and positive with $\lambda'(0)>0$, then 
the formulae  (\ref{mg3},\ref{mbg3}) define  Riemannian metrics. 
\end{Rem}

Proof:  Let the metrics be  proportional at 
$x_0$. Take a small $\epsilon-$ball $B_\epsilon$ (in the metric $g$) 
with the center in $x_0$.  If $\epsilon$ is small enough, $x_0$ is the 
only point of the ball where the metrics are proportional.  Indeed, suppose   
they are proportional at three points $y_1,y_2,y_3$ of the ball. 
Then, for almost every  point $x$ 
 of the ball, there exist three geodesics
 $\gamma_1,\gamma_2,\gamma_3$ such that 
$\gamma_i(0)=y_i$, $\gamma_i(1)=x$ and the velocity vectors $\dot\gamma_1(1),\dot\gamma_2(1),\dot\gamma_3(1)$ are mutually transversal. Let us show that  
$\lambda_1=\lambda_2$ or $\lambda_2=\lambda_3$ at $x$.
 Indeed, by Lemma~\ref{technical}, for every  $i=1,2,3$ we have that 
 $\lambda_2$ 
is a double root of 
$I_t(\gamma_i(1),\dot\gamma_i(1))$. 
At $T_xB_\epsilon$, consider the coordinate system where $g$ and $L$
are given by diagonal matrices $\mbox{diag}(1,1,1)$ and 
$\mbox{diag}(\lambda_1,\lambda_2,\lambda_3)$, respectively.
In this coordinate system, the polynomial $I_t$ is given by 
\begin{equation} \label{p3} 
(\lambda_2-t)(\lambda_3-t)\xi_1^2+ (\lambda_1-t)(\lambda_3-t)\xi_2^2+(\lambda_1-t)(\lambda_2-t)\xi_3^2. 
\end{equation}
Then, the components  $(\xi_1,\xi_2,\xi_3)$ 
 of the velocity vectors $\dot\gamma_i(1)$ satisfy 
\begin{equation}\label{36}
\left\{
\begin{array}{ccc}
(\lambda_1(x)-\lambda_2)(\lambda_3(x)-\lambda_2)\xi_2^2 &=& 0\\
(\lambda_1(x)-\lambda_2)(\xi_2^2+\xi_3^2)+
(\lambda_3(x)-\lambda_2)(\xi_2^2+\xi_1^2)&=&0.  
\end{array}
\right. 
\end{equation} 
If $\lambda_1(x)<\lambda_2<\lambda_3(x)$, the solutions of 
(\ref{36}) are organized into two intersected straight lines so that 
the  velocity vectors $\dot\gamma_i(1)$ are not mutually transversal. 
Then, at almost every  point of $B_\epsilon$ we have 
$\lambda_1=\lambda_2$ or $\lambda_2=\lambda_3$, 
which contradicts Corollary~\ref{ordered12}. 

Consider the set 
$$
U\eqdef \{ x\in B_\epsilon: \ (\lambda_1(x)-\lambda_2)(\lambda_{3}(x)-\lambda_2)=0\}.
$$

\begin{Lemma}\hspace{-2mm}{\bf .} 
The set $U$ is a totally geodesic connected submanifold of 
$B_\epsilon$  of dimension 1. (In other words, $U$ is a geodesic segment). 
\end{Lemma} 
{Proof:} First of all, there exists $x_1\in B_\epsilon$ where 
precisely two eigenvalues of $L$ coincide. Indeed, take the 
$\frac{\epsilon}{2}-$ sphere $S_{\frac{\epsilon}{2}}$ 
with the center in $x_0$ and consider the exponential 
mapping $\mbox{exp}:T_{x_0}B_\epsilon\to B_\epsilon$. Suppose there 
exists no point of  $S_{\frac{\epsilon}{2}}$   where 
precisely two eigenvalues of $L$ coincide. Then, at every point of 
the sphere $S_{\frac{\epsilon}{2}}$ the eigenspace of $L$ corresponding to 
$\lambda_2$ has dimension one. Let us show that it is tangent to the sphere. 
It is sufficient to show that it is orthogonal to the geodesic 
connecting the point of the sphere with the point $x_0$. Denote the 
initial velocity vector of the geodesic by $\xi$.  
At the tangent space at this point, 
 consider the coordinate system where $g$ and $L$
are given by diagonal matrices $\mbox{diag}(1,1,1)$ and 
$\mbox{diag}(\lambda_1,\lambda_2,\lambda_3)$, respectively.
In this coordinate system, the polynomial $I_t$ is given by (\ref{p3}).
Since the geodesic goes through the point where 
$\lambda_1=\lambda_2=\lambda_3$, by Lemma~\ref{technical},  $\lambda_2$ is 
the double-root $I_t(\xi)$. Then the components $\xi_1,\xi_2,\xi_3$ of 
$\xi$ satisfy the system (\ref{36}). From the first equation of the system, 
in  view that $\lambda_1<\lambda_2<\lambda_3$, we obtain that the component $\xi_2$ is zero so that $\xi$ is orthogonal to the eigenspace of $L$ corresponding to $\lambda_2$. Thus, the eigenspaces of $L$ corresponding to $\lambda_2$ 
give us a smooth one-dimensional distribution on the 2-sphere which is
 impossible because the Euler characteristic of the sphere. Finally, 
there exists  $x_1\in S_{\frac{\epsilon}{2}}$ where precisely two
 eigenvalues of $L$ coincide.

 Denote by $\gamma$ the  geodesic going through  $x_1$ and $x_0$.
 Let us show that at every point of this geodesic at least  
 two
 eigenvalues of $L$ coincide.

We  assume that  $\gamma(1)=x_1$ and $\gamma(0)=x_0$.

 At each tangent space, we 
   can find coordinates such that  
  $g$ and  $L$  are 
 given by the  diagonal matrixes $\diag(1,1,1)$ and 
$\diag(\lambda_1,\lambda_2,\lambda_3)$, respectively.  
Consider the function 
$$
I'_t \eqdef -\left(\frac{d}{dt}I_t\right).
$$
For every  fixed $t$, the 
function $I'_t$  is  an integral of the geodesic flow. 
By (\ref{p3}),  
\begin{eqnarray}
I'_t &\hspace{-4mm} =  & \hspace{-4mm} (\lambda_2-t+\lambda_3(x)-t)\xi_1^2+(\lambda_1(x)-t+\lambda_3(x)-t)\xi_2^2
+  (\lambda_2-t+\lambda_1(x)-t)\xi_3^2, \nonumber \\ 
 I'_{\lambda_2}&\hspace{-4mm}=&  \hspace{-4mm}
(\lambda_1(x)-\lambda_2)(\xi_2^2+\xi_3^2)+
(\lambda_3(x)-\lambda_2)(\xi_2^2+\xi_1^2). \label{34}
\end{eqnarray}
 
By Lemma~\ref{technical},  
   $\lambda_2$ is a 
double-root of the polynomial $I_t(\gamma,\dot\gamma)$. It follows from 
 Lemma~\ref{technical} that   the 
leading coefficient of the polynomial $I'_{t}(\gamma,\dot\gamma)$
 is not zero. 

Let us prove that the  
differential $dI'_{\lambda_2}$ 
vanishes   at each point of the 
 geodesic orbit $(\gamma,\dot\gamma)$. Since $I'_{\lambda_2}$
is an integral, it is sufficient to show this at the point  
$(\gamma(1),\dot\gamma(1))$ only. 
We have,  
$$
I'_{\lambda_1(x)}=(\lambda_3(x)-\lambda_1(x))(\xi_1^2+\xi_2^2)+
(\lambda_2-\lambda_1(x))(\xi_1^2+\xi_3^2). 
$$ 
We see that the function 
$I'_{\lambda_1(x)}$
is non-negative. 
At the point $(\gamma(1),\dot\gamma(1))$, it coincides with $I'_{\lambda_2}$
and, therefore,  is zero. 
Then,  it 
has a minimum at the point $(\gamma(1),\dot\gamma(1))$, and its 
differential vanishes. 

Let us show that the differential 
of the function $I'_{\lambda_1(x)}-I'_{\lambda_2}$  also  
vanishes  at the point 
$(\gamma(1),\dot\gamma(1))$. 
Indeed, the function  $I'_t$ is a linear polynomial  in
 $t$  with non-zero leading coefficient at the point 
$(\gamma(1),\dot\gamma(1))$. 
  Since $\lambda_1(x)<\lambda_2$, then 
 the function 
 $I'_{\lambda_1(x)}-I'_{\lambda_2}$ is either everywhere positive or 
everywhere negative. Since it vanishes  at  
$(\gamma(1),\dot\gamma(1))$,
 the differential of   
$I'_{\lambda_1(x)}-I'_{\lambda_2}$ vanishes at the
 point 
$(\gamma(1),\dot\gamma(1))$. 

Thus  the differential 
 $dI'_{\lambda_2}$ 
is zero at each point of the geodesic orbit $(\gamma,\dot\gamma)$. 
At each point of $\gamma$, the
components  $\frac{\partial I'_{\lambda_2}}{\partial \xi_i}$     of 
  $dI'_{\lambda_2}$ are  
$$
2(\lambda_3-\lambda_2)\xi_1, \ 
2(\lambda_3-\lambda_2-\lambda_2+\lambda_1)\xi_2,  \ 
2(\lambda_1-\lambda_2)\xi_3.
$$
Since the differential vanishes, all its components  are equal to zero. 
Then,   $\lambda_2=\lambda_3$ or 
$\lambda_1=\lambda_2$ or $\xi_2\ne 0$. 

On the other hand, by (\ref{p3}), using that $I_{\lambda_2}=0$, we see that  
$ (\lambda_2-\lambda_1(x))(\lambda_2-\lambda_3(x))=0$ or $\xi_2=0$. 
  
Finally, every   point  of $\gamma$ 
lies in $U$. 

Now let us prove that $\gamma$ actually coincides with $U$. Assume 
$\gamma$ does not coincide with $U$; that is, there exist a point $x_2$
where $\lambda_1=\lambda_2$ or $\lambda_2=\lambda_3$ not lying on
 the geodesic $\gamma$. Then, for almost every  $x\in B_\epsilon$,  
there exist three geodesics $\gamma_1,\gamma_2,\gamma_3$ such that 
$\gamma_i(1)=x$ and $\gamma_i(0)\in U$ and the vectors 
$\dot\gamma_i(1)=x$ are linearly independent. From other side, 
the components
 $\xi_1,\xi_2,\xi_3$ of each of these vectors (in the coordinate system 
 where $g$ and $L$
are given by diagonal matrices $\mbox{diag}(1,1,1)$ and 
$\mbox{diag}(\lambda_1,\lambda_2,\lambda_3)$) satisfy the first equation of 
the system (\ref{36}). Then, $x\in U$ which contradicts 
Lemma~\ref{technical}. 
Finally, $U$ coincides with $\gamma$.

Now let us show that there exists a smooth 
vector field $v_2$ on $B_\epsilon$ 
such that: 
\begin{enumerate} 
\item $v_2$ is a Killing vector field with respect to both metrics, 
\item $g(v_2,v_2)=(\lambda_1-\lambda_2)(\lambda_1-\lambda_3)$  \label{2}
(in particular, $v_2$ vanishes at $U$),  
\item $Lv_2=\lambda v_2$, \label{3} 
\item the integral curves of $v_2$ are homeomorphic to the circle. 
\end{enumerate} 

Let us prove that there exists $v_2$ on 
  $B_\epsilon\setminus U$  satisfying (\ref{2}),(\ref{3}).  Since at every  point of $B_\epsilon\setminus U$  the eigenspace of $L$ corresponding to $\lambda_2$ 
is one-dimensional, at every  point there exist precisely two 
vectors  satisfying conditions  (\ref{2}),(\ref{3}).

Take a point $x_2\in B_\epsilon\setminus U$ and choose a vector  satisfying conditions  (\ref{2}),(\ref{3}) at $x_2$. 
 Then  we can smoothly   
choose the vector field  satisfying (\ref{2}),(\ref{3}) locally, and, therefore, along any curve. Thus it is sufficient to show that the result does not depend on the curve. Take two curves $c_1$, $c_2$
 connecting $x_2$ and an arbitrary point $x \in  B_\epsilon\setminus U$.   

Take a point $x_1\in U$, $x_1\ne x_0$. 
Without loss of generality we can assume that $\lambda_1(x_1)<\lambda_2(x_1)=\lambda_3(x_1)$. By Corollary~\ref{l3} 
in a small neighborhood $W$  
of $x_1$ the metrics are given by the model form  
(\ref{l4},\ref{l5}). 
Clearly, 
there exists curves $\tilde c_1$, 
$\tilde c_2$ such that $c_1$ is homotop to $\tilde c_1$, 
$c_2$ is homotop to $\tilde c_2$, 
 and $\tilde c_1$ coincides with   
$\tilde c_2$ outside of $W$. As it follows from Corollary~\ref{l3}, 
we can choose $v_2$ satisfying (\ref{2}),(\ref{3}) in $W$.
Since the curves coincide outside of $W$, there exists $v_2$ 
on $B_\epsilon\setminus U$   satisfying (\ref{2}),(\ref{3}).   

In Levi-Civita coordinates from Theorem~\ref{Levi-Civita}, 
the vector field $v_2$ equals $\pm \frac{\partial }{\partial x_2}$ 
and therefore is Killing with respect to the metrics $g$ and $\bar g$.  
Put the vector field $v_2$ equal to zero at every point of $U$. 
Since any isometry of $g$ is a diffeomorphism, the vector field is smooth 
everywhere. Consider the exponential mapping from the point $x_0$. 
As we have proven before, the  $v_2$ is tangent to the images of the 
spheres on $T_{x_0}B_\epsilon$, and, therefore, generates a Killing 
vector field on every such sphere. Then its integral curves are closed 
\cite{BMF}.   

Consider the point $x_1$ where  $\lambda_1(x_1)<\lambda_2(x_1)=
\lambda_3(x_1)$.  In the neighborhood of the point there exists coordinates 
$x_1,x_2,x_3$ where the metrics are given by the form (\ref{l4},\ref{l5}).
Consider the exponential mapping $exp:T_{x_1}B_\epsilon\to B_\epsilon$, the 
2-plane $P\subset T_{x_1}B_\epsilon$
 spanned by $ \frac{\partial }{\partial x_1}$ and 
$\frac{\partial }{\partial x_3}$, and the image $exp(P)$ 
 of this plane under the 
 exponential mapping. 

Since $\epsilon$ is small, $exp(P)$ is a two-dimensional submanifold 
of $B_\epsilon$. 
Let us show that  
$exp(P)$  is  totally geodesic,
 (so that the restrictions of the metrics to $exp(P)$
  are projectively equivalent), that the eigenvalues of
 the tensor $L$ constructed for the restriction of the metrics 
to $exp(P)$ are $\lambda_1,\lambda_3$, and that at every point  $exp(P)$
is orthogonal to the vector $v_2$. 

It is sufficient to prove both facts at almost every point of  $exp(P)$. 
Take a point $p_2\in exp(P)$ where $\lambda_1<\lambda_2<
\lambda_3$. Denote by $\gamma$ the geodesic connecting this point with $p_1$, 
$\gamma(1)=p_1$,  $\gamma(2)=p_2$. In a neighborhood of $p_1$, the 
listed above statements are true because of Corollary~\ref{l3}. In a neighborhood of another points,  the 
listed above statements are true because of Levi-Civita's Theorem. 

Finally, since  $exp(P)$ is totally geodesic, the restrictions of the 
metrics $g$ and $\bar g$ to $exp(P)$ are projectively equivalent. They are 
proportional at $x_0$, they are 
 non-proportional at  every other point, and the eigenvalues of $L$ constructed for the restriction are non-constant. By Theorem~\ref{local2} the metrics 
are  given by (\ref{l0},\ref{l1}). Finally, the metrics $g$ and $\bar g$ 
 are as in Theorem~\ref{local3}.

\subsection{ General case   and realization} \label{2.6}
Let $M^{n_1}_1$, $  M^{n_2}_2$ be connected  manifolds. 
 Suppose Riemannian metrics $g_1,\bar g_1$  on $M^{n_1}_1$ are   projectively equivalent
and are  strictly-non-proportional at least at one point.
Suppose Riemannian metrics $g_2,\bar g_2$  on $M^{n_2}_2$ are  
 projectively equivalent
and are  strictly-non-proportional at least at one point.  We denote by $L_1$ the tensor (\ref{l}) corresponding to $g_1,\bar g_1$, and by  $L_2$ the tensor (\ref{l}) corresponding to 
$g_2,\bar g_2$. Assume in addition 
  the eigenvalues of $L_1$ are less than the eigenvalues   of $L_2$. 
Consider the direct product $M^{n_1}_1\times M^{n_2}_2$
with the canonical transversal foliations $B_{n_1}$ and 
$B_{n_2}$. (The leaves of the foliation $B_{n_1}$ are the products of 
$M_1^{n_1}$ and a point of $M_2^{n_2}$,  the 
leaves of the foliation $B_{n_2}$ are the products of the  point of $M_1^{n^1}$
and $M_2^{n_2}$.)
Consider the Riemannian metric $h\eqdef g_1+g_2$ and 
$\bar h\eqdef \bar g_1+\bar g_2$ on the  product $M^{n_1}_1\times M^{n_2}_2$.
 It is easy to see that the foliations and the metrics satisfy the assumptions 
(i,ii,iii) of Section~\ref{2.3}. Then, by Lemma~\ref{inverse}, 
the metrics $g$ and $\bar g$ given 
by (\ref{g},\ref{bg}) are projectively equivalent and are 
strictly-non-proportional at least at one point of the manifold.  
Thus, given  two triples $( M^{n_1}_1, g_1,\bar g_1)$, 
$( M^{n_2}_2, g_2,\bar g_2)$, we constructed the triple $(M^{n_1}_1\times M^{n_2}_2, g, \bar g)$. We will denote this operation by ``$\oplus$'': 
$$
(M^{n_1}_1, g_1,\bar g_1)\oplus 
( M^{n_2}_2, g_2,\bar g_2)=(M^{n_1}_1\times M^{n_2}_2, g, \bar g).
$$
It is easy to check that the operation is associative: 
$$
((M^{n_1}_1, g_1,\bar g_1)\oplus 
( M^{n_2}_2, g_2,\bar g_2))\oplus (M^{n_3}_3, g_3,\bar g_3)=(M^{n_1}_1, g_1,\bar g_1)\oplus 
(( M^{n_2}_2, g_2,\bar g_2)\oplus (M^{n_3}_3, g_3,\bar g_3)).
$$

Consider connected  manifolds $M_1^{k_1}, M_2^{k_2},...,M_m^{k_m}$ 
with   projectively equivalent
 metric $g_1,\bar g_1;g_2,\bar g_2;...;g_m,\bar g_m$, respectively. 
Assume the metrics  are strictly-non-proportional at least at one point. 
We denote by $L_i$ the tensor (\ref{l}) constructed for $g_i,\bar g_i$, 
$i=1,...,m$. 
Assume in addition that,  for any $i<j$, 
  the eigenvalues of $L_i$ are less than the eigenvalues   of $L_j$. 
 Then,  we can canonically construct projectively equivalent
 metric on the product of the manifolds, and these metrics are 
 strictly-non-proportional at least at one point.

\begin{Cor}\hspace{-2mm}{\bf .} \label{product} 
A product of spheres admits projectively equivalent metrics
which are strictly-non-proportional at least at one point.  
\end{Cor}

\noindent {\bf Proof:}  Basically we will show that if 
connected   closed manifolds $M_1^{k_1}, M_2^{k_2},...,M_m^{k_m}$ 
admit   projectively equivalent strictly-non-proportional at least at one 
point 
 metric $g_1,\bar g_1;g_2,\bar g_2;...;g_m,\bar g_m$, then the product of the 
$M_1^{k_1}\times  M_2^{k_2}\times ... \times M_m^{k_m}$ also admit  projectively equivalent metrics strictly-non-proportional at least at one 
point.  
  
Since the  manifolds $M_i^{k_i}$  are closed,  the 
eigenvalues of the tensors (\ref{l}) constructed for the metrics 
$g_i, \bar g_i$ are bounded. By definition,  
the tensor  (\ref{l}) constructed for  the metrics  $g_i, \bar g_i$ and 
the tensor  (\ref{l}) constructed for  the metrics  $C\cdot g_i, \bar g_i$,
 where 
$C$ is a positive constant, are related by
$$
L_{new}=C^{\frac{1}{{k_i}+1}}L_{old}.
$$
Thus  without loss of generality we can 
assume that for  $i<j$ 
the eigenvalues of  the tensor  (\ref{l}) 
constructed for  the metrics  $g_j, \bar g_j$
are greater than the eigenvalues of 
the tensor  (\ref{l}) constructed for  the metrics  $g_i, \bar g_i$. 
Then,  by Lemma~\ref{inverse}, 
$$
(M^{k_1}_1, g_1,\bar g_1)\oplus 
(M^{k_2}_2, g_2,\bar g_2))\oplus...\oplus (M^{k_m}_3, g_m,\bar g_m)
$$ gives 
us projectively equivalent metrics on $ M_1^{k_1}\times  M_2^{k_2}\times ... \times M_m^{k_m}$  strictly-non-proportional at least at one 
point.

Finally, in order to prove Corollary~\ref{product}, we need to show that 
the sphere $S^n$ admits two projectively equivalent metrics which are 
strictly-non-proportional at least at one point. Essentially, it was done in 
\cite{Beltrami}: 
 the metric 
$g$ is the restriction of the Euclidean metrics $dx_1^2+...+dx_{n+1}^2$ 
to the sphere 
$$
S^n\eqdef
 \{(x_1,x_2, ...,x_{n+1})\in R^{n+1}: \ x_1^2+x_2^2+...+x_{n+1}^2=1\}. 
$$
The metric $\bar g$ is the pull-back 
$l^*g$, where the mapping $l:S^n\to S^n$ is given by 
$l:v\mapsto \frac{A(v)}{\|A(v)\|}$, where $A$ is 
an arbitrary linear non-degenerate transformation of $R^{n+1}$. 

The metrics $g$ and $\bar g$ are projectively equivalent. 
Indeed, the geodesics of  $g$ are great circles 
(the intersections of planes that go through the origin with the 
 sphere). The mapping $A$ is linear and, hence,
 takes  
planes to  planes. Since the normalization 
$w\mapsto \frac{w}{\|w\|}$ takes   planes to their intersections 
with the sphere, the mapping $l$ takes  great circles to  
great circles. Thus  the metrics $g$ and $\bar g$
are projectively equivalent. It is easy  to verify that 
for  almost all  linear transformations  $A$,   
(and in particular 
 for $A=\diag(a_1,...,a_{n+1})$, where $a_1<a_2<....<a_{n+1}$), 
 the metrics   $g$ and $\bar g$ are strictly-non-proportional at
 almost  every point.  Corollary~\ref{product} is proved. 

\begin{Def}\hspace{-2mm}{\bf .} {\rm 
By {\it prime   standard  triples}  we will  mean:
\begin{itemize} 
\item In the one-dimensional case, 
 $(I, dx^2,\frac{1}{\lambda(x)}dx^2)$, where $I$ is an interval with the 
coordinate $x$ and $\lambda$ is a smooth positive function. 

\item In the two-dimensional case, $(D^2, g, \bar g)$, where  $D^2$ is a 2-disc
with coordinates $u,v$ and the metrics $g, \bar g$ are given either by the 
formulae (\ref{mg},\ref{mbg}), (\ref{l0},\ref{l1}) or (\ref{l0},\ref{l2}), 
respectively. 

\item In the three-dimensional case, $(D^3, g, \bar g)$, where  $D^3$ is a 3-disc with coordinates $u_1,u_2,u_3$
 and the metrics $g, \bar g$ are given  by the 
formulae (\ref{l4},\ref{l5}). 
\end{itemize} 

A  disk $D^n$ with two projectively equivalent metrics 
$g$, $\bar g$ 
 will be called {\it standard} is there exists the 
prime standard triples $(D_i,g_i,\bar g_i)$, $i=1,...,m$, 
 such that 
$$
(D_1, g_1,\bar g_1)\oplus 
(D_2, g_2,\bar g_2))\oplus...\oplus (D_m, g_m,\bar g_m)=(D^n, g,\bar g). 
$$
The sign ``='' here mean the existence of a diffeomorphism between
 $D_1\times...\times D_m$ and $D^n$ which is an isometry with respect 
to both metrics. }
\end{Def}

\begin{Th}\hspace{-2mm}{\bf .} \label{local}
Suppose projectively equivalent  Riemannian metrics $g$ and $\bar g$ on a connected manifold $M^n$ 
and are strictly non-proportional at least at one point. Then, every 
 point 
$p$ of the manifold has a neighborhood $U$ such that the 
triple $(U, g, \bar g)$  is standard. 
\end{Th}

\noindent{\bf Proof:} We use induction in dimension. 
 For dimension  one the statement is trivial; for dimension two, 
 the statement  follows from Theorem~\ref{local2}. 
Assume $n\ge 3$, and  
suppose the  theorem is true for dimensions less than $n$. 
Let us prove that it is true for dimension $n$.  By Corollary~\ref{ordered12}, 
 there exists $1 \le r \le 3$ such that 
$\lambda_1(p)=...=\lambda_r(p)<\lambda_{r+1}(o)\le ...\le\lambda_n(p)$. If $r=n$, then  $n=3$; then, 
 the theorem follows from Theorem~\ref{local3}.

Suppose $n>r$. Then, it is so at every  point of a small neighborhood of $p$, 
and we can apply the splitting procedure from Section~\ref{2.2}. We obtain the
metrics $h,\bar h$ and the foliations $B_r$ and $B_{n-r}$ such that 
the local product structures 
$(h, B_r,   B_{n-r})$ and $(\bar h, B_r,   B_{n-r})$ satisfy the conditions 
(i,ii,iii) of Section~\ref{2.3}. By definition of a local product structure, 
  a (possibly,  smaller) neighborhood of $p$, is
a direct product of two discs $D^r$ and $D^{n-r}$,  the metrics $h$ and 
$\bar h$ are the product metrics $g_1+g_2$, $\bar g_2+\bar g_2$, where 
$g_i$ is projectively equivalent to $\bar g_i$. 
By the induction assumption, the triples $(D^r, g_1,\bar g_1)$ and  
$(D^{n-r}, g_2,\bar g_2)$ are standard; therefore, 
$(D^r, g_1,\bar g_1)\oplus (D^{n-r}, g_2,\bar g_2)$ is standard as well. 
In view of Remark~\ref{missleading},  $(D^r, g_1,\bar g_1)\oplus (D^{n-r}, g_2,\bar g_2)$ is precisely $(D^r\times D^{n-r}, g, \bar g)$. Theorem \ref{local} 
 is proved.

\end{document}